\newif\ifpdf
\ifx\pdfoutput\undefined
        \pdffalse
\else
        \pdfoutput=1 
        \pdftrue
\fi

\ifpdf
\documentclass[pdftex]{amsart}
\else
\documentclass{amsart}
\fi

\usepackage{amssymb}

\ifpdf
        \usepackage{color}
        \usepackage[pdftex]{graphicx}
\else
        \usepackage{graphicx}
\fi

\usepackage[all]{xy}

\title[Parall{\'e}lisme et homotopie]{Concurrent Process up to Homotopy (I)}
\author[P. Gaucher]{Philippe Gaucher}
\address{Institut de Recherche Math\'ematique Avanc\'ee\\ ULP et
CNRS\\ 7 rue Ren\'e Descar\-tes\\ 67084 Strasbourg Cedex\\ France}
\email{gaucher@math.u-strasbg.fr}
\urladdr{http://www-irma.u-strasbg.fr/\~{}gaucher/}

\newcommand{\p}\times
\renewcommand{\vec}{\overrightarrow}
\renewcommand{\P}{\mathbb{P}}

\newcommand{\be}{\begin{equation}}
\newcommand{\ee}{\end{equation}}
\newcommand{\bea}{\begin{eqnarray}}
\newcommand{\eea}{\end{eqnarray}}
\newcommand{\beas}{\begin{eqnarray*}}
\newcommand{\eeas}{\end{eqnarray*}}

\newtheorem{theoreme}{Th{\'e}or{\`e}me}[section]
\newtheorem{prop}[theoreme]{Proposition}

\newtheorem{rem}[theoreme]{Remarque}

\newtheorem{definition}[theoreme]{D{\'e}finition}

\newcommand{\bd}{\begin{defn}}
\newcommand{\ed}{\end{defn}}
\newcommand{\bcd}{\begin{defn}}
\newcommand{\ecd}{\end{defn}}
\newcommand{\bex}{\begin{exmp}}
\newcommand{\eex}{\end{exmp}}
\newcommand{\bp}{\begin{prop}}
\newcommand{\ep}{\end{prop}}
\newcommand{\bth}{\begin{theoreme}}
\renewcommand{\eth}{\end{theoreme}}
\newcommand{\br}{\begin{rem}}
\newcommand{\er}{\end{rem}}
\newcommand{\bpf}{\begin{proof}}
\newcommand{\epf}{\end{proof}}

\newcommand{\fl}[1]{\ar@{->}[l]_{#1}}
\newcommand{\fr}[1]{\ar@{->}[r]^-{#1}}
\newcommand{\fd}[1]{\ar@{->}[d]_{#1}}
\newcommand{\fu}[1]{\ar@{->}[u]^{#1}}
\newcommand{\f}[2]{\ar@{->}[#1]|{#2}}
\newcommand{\ff}[2]{\ar@2{->}[#1]|{#2}}
\newcommand{\frr}[1]{\ar@{->}[rr]^{#1}}

\newcommand{\ho}{{\mathbf{Ho}}}

\newcommand{\iso}{\cong}

\newcommand{\vI}{\vec{I}}

\renewcommand{\leq}{\leqslant}
\renewcommand{\geq}{\geqslant}

\def\cartesien{%
  \ar@{-}[]+R+<6pt,-2pt>;[]+RD+<6pt,-6pt>%
  \ar@{-}[]+D+<2pt,-6pt>;[]+RD+<6pt,-6pt>%
}
\def\cocartesien{%
  \ar@{-}[]+L+<-6pt,+2pt>;[]+LU+<-6pt,+6pt>%
  \ar@{-}[]+U+<-2pt,+6pt>;[]+LU+<-6pt,+6pt>%
}

\newcommand{\brm}[1]{\rm{\mathbf{#1}}}

\newcommand{\gltop}{{\brm{glTop}}}

\newcommand{\dtop}{{\brm{Flow}}}

\newcommand{\tdtop}{{\brm{FLOW}}}
\newcommand{\ttop}{{\brm{TOP}}}
\newcommand{\tgltop}{{\brm{glTOP}}}

\newcommand{\limind}{\varinjlim}
\newcommand{\limproj}{\varprojlim}

\makeatletter

\def\varholim@#1#2{%
  \vtop{\m@th\ialign{##\cr
    \hfil$#1\operator@font holim$\hfil\cr
    \noalign{\nointerlineskip\kern1.5\ex@}#2\cr
    \noalign{\nointerlineskip\kern-\ex@}\cr}}%
}
\def\holimproj{%
  \mathop{\mathpalette\varholim@{\leftarrowfill@\textstyle}}\nmlimits@
}
\def\holimind{%
  \mathop{\mathpalette\varholim@{\rightarrowfill@\textstyle}}\nmlimits@
}

\newskip\@bigflushglue \@bigflushglue = -100pt plus 1fil

\def\bigcentering{\let\\\@centercr\rightskip\@bigflushglue%
\leftskip\@bigflushglue
\parindent\z@\parfillskip\z@skip}

\makeatother

\DeclareMathOperator{\ctop}{{\textbf{Top}}}

\begin{document}

\begin{abstract}
Les CW-complexes globulaires et les flots sont deux mod{\'e}lisations g{\'e}om{\'e}triques des
automates parall{\`e}les qui permettent de formaliser la notion de dihomotopie. La dihomotopie
est une relation d'{\'e}quivalence sur les automates parall{\`e}les qui pr{\'e}serve des propri{\'e}t{\'e}s
informatiques comme la pr{\'e}sence ou non de deadlock. On construit un plongement des
CW-complexes globulaires dans les flots et on d{\'e}montre que deux CW-complexes
globulaires sont dihomotopes si et seulement si les flots associ{\'e}s sont dihomotopes.

Globular CW-complexes and flows are both geometric models of
concurrent processes which allow to model in a precise way the notion
of dihomotopy. Dihomotopy is an equivalence relation which preserves
computer-scientific properties like the presence or not of
deadlock. One constructs an embedding from globular CW-complexes to
flows and one proves that two globular CW-complexes are dihomotopic if
and only if the corresponding flows are dihomotopic.
\end{abstract}
\maketitle

\section{Rappels sur les CW-complexes globulaires}

Cette note est la premi{\`e}re de deux notes pr{\'e}sentant quelques r{\'e}sultats de
\cite{flow}. Tous les espaces topologiques sont suppos{\'e}s faiblement 
s{\'e}par{\'e}s et compactement engendr{\'e}s, c'est-{\`a}-dire dans ce cas
hom{\'e}omorphes {\`a} la limite inductive de leurs sous-espaces compacts
(cf. l'appendice de \cite{Ref_wH} pour un survol des propri{\'e}t{\'e}s de ces
espaces). On travaille ainsi dans une cat{\'e}gorie d'espaces topologiques
(not{\'e}e $\ctop$) qui est non seulement compl{\`e}te et cocompl{\`e}te mais en
plus cart{\'e}siennement ferm{\'e}e \cite{MR35:970}. En d'autres termes, le
foncteur $-\p X:\ctop\rightarrow \ctop$ a un adjoint {\`a} droite
$\ttop(X,-)$.  Dans la suite, pour $n\geq 1$, $D^n$ est le disque
ferm{\'e} de dimension $n$ et $S^{n-1}$ est le bord de $D^n$, {\`a} savoir la
sph{\`e}re de dimension $n-1$. En particulier la sph{\`e}re de dimension $0$
est la paire $\{-1,+1\}$.

Si $Z$ est un espace topologique non vide, on note $Glob^{top}(Z)$
l'espace topologique obtenu en partant de $Z\p [0,1]$ et en
quotientant par les relations $(z,0)=(z',0)$ et $(z,1)=(z',1)$ pour
tout $z,z'\in Z$. Cet espace est muni de l'ordre partiel suivant :
$(z,t)\leq (z',t')$ si et seulement si $z=z'$ et $t\leq t'$. La classe
d'{\'e}quivalence de $(z,0)$ (resp. $(z,1)$) est inf{\'e}rieure
(resp. sup{\'e}rieure) {\`a} tous les points de $Glob^{top}(Z)$.  Le
\textit{$0$-squelette} $Glob^{top}(Z)^0$ de $Glob^{top}(Z)$ sera la paire constitu{\'e}e de la
classe de $(z,0)$ et de celle de $(z,1)$. On pose $\vI=Glob^{top}(\{*\})$
qui a donc pour $0$-squelette $\{(*,0),(*,1)\}$.

Un \textit{CW-complexe globulaire} $X$ est, par d{\'e}finition, obtenu
comme limite inductive $\limind_n X^1_n$ d'espaces topologiques de la
fa{\c c}on suivante. On part d'un espace discret $X^0$ appel{\'e} le
\textit{$0$-squelette} et on commence par lui attacher des $\vI$ pour
obtenir un premier espace topologique $X^1_0$. Puis on suppose $X^1_n$
construit pour $n\geq 0$. On obtient alors $X^1_{n+1}$ en attachant {\`a}
$X^1_n$ des $Glob^{top}(D^{n+1})$ le long de $Glob^{top}(S^{n})$. Tout cela de
telle fa{\c c}on que les attachements des $\vI$ et les morphismes
d'attachement $Glob^{top}(S^{n})\rightarrow X^1_n$ soient des morphismes
de CW-complexes globulaires (cf. ci-dessous).

\begin{figure}
\ifpdf \centerline{\pdfimage width 5cm {exglob.png}} \else
\begin{center}
\includegraphics[width=5cm]{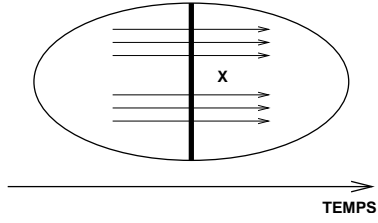}
\end{center}
\fi
\caption{Repr{\'e}sentation symbolique de $Glob^{top}\left(X\right)$ pour un espace topologique $X$} \label{exglob}
\end{figure}

Un \textit{chemin d'ex{\'e}cution} $\gamma:\vI\rightarrow X$ est une
application continue localement croissante pour l'ordre partiel d{\'e}fini
sur chaque globe de $X$ et telle que $\gamma(\vI^0)\subset
X^0$. L'ensemble de ces $\gamma$ \textit{non-constants} est not{\'e}
$\mathbb{P}^{top}(X)$. Un morphisme de CW-complexes globulaires
$f:X\rightarrow Y$ pour $X$ et $Y$ quelconques est une application
continue telle que $f(X^0)\subset Y^0$ et telle que $\gamma\mapsto
f\circ \gamma$ induise une application de $\mathbb{P}^{top}(X)$ dans
$\mathbb{P}^{top}(Y)$. L'hypoth{\`e}se de non-contractibilit{\'e} des chemins
d'ex{\'e}cution par $f$, qui peut para{\^i}tre {\'e}trange en premi{\`e}re lecture,
est en fait essentielle pour pouvoir faire certaines constructions
homologiques : cf. par exemple le chapitre "why non-contracting maps
?" de \cite{diCW} ou bien \cite{Gau}.

\begin{definition}\cite{diCW}
Deux morphismes de CW-complexes globulaires $f$ et $g$ de $X$ dans $Y$ sont S-homotopes s'il existe
une application continue $H:X\p [0,1]\rightarrow Y$ telle que
\begin{enumerate}
\item $H(-,u)$ est un morphisme de CW-complexes globulaires de $X$ dans $Y$ pour tout $u\in [0,1]$
\item $H(-,0)=f$ et $H(-,1)=g$
\end{enumerate}
On {\'e}crit $f\sim_S g$.
\end{definition}

En particulier deux morphismes de CW-complexes globulaires
S-homotopes co{\"\i}n\-cident sur le $0$-sque\-let\-te.  De l{\`a} on d{\'e}finit
la notion de CW-complexes globulaires S-homotopes :

\begin{definition}\cite{diCW}
Deux CW-complexes globulaires $X$ et $Y$ sont S-homotopes s'il
existe un morphisme de CW-complexes globulaires $f:X\rightarrow Y$ et
un morphisme de CW-complexes globulaires $g:Y\rightarrow X$ telle que
$f\circ g\sim_{S} Id_Y$ et $g\circ f\sim_{S} Id_X$. \end{definition}

\begin{definition}\cite{diCW}
Une T-homotopie $f:X\rightarrow Y$ est un morphisme de CW-complexes
globulaires qui induit un hom{\'e}omorphisme entre les espaces topologiques
sous-jacents. On dit alors que $X$ et $Y$ sont T-homotopes. \end{definition}

On peut alors d{\'e}montrer le

\begin{theoreme}\cite{diCW}
La localisation $\ho(\gltop)$ de la cat{\'e}gorie $\gltop$ des CW-complexes globulaires par
rapport aux S-homotopies et aux T-homotopies existe. \end{theoreme}

La forme de globe
suffit pour mod{\'e}liser les automates parall{\`e}les \cite{Pratt}. De plus la
S-homotopie et la T-homotopie ne modifient pas des propri{\'e}t{\'e}s
informatiques comme la pr{\'e}sence ou non de deadlock \cite{diCW}. La cat{\'e}gorie $\ho(\gltop)$ est
donc un cadre possible pour l'{\'e}tude des automates parall{\`e}les {\`a} homotopie pr{\`e}s.

\section{Flot}

\begin{definition}\cite{flow} Un flot $X$ consiste en la donn{\'e}e d'un ensemble $X^0$ appel{\'e} $0$-squelette,
d'un espace topologique $\P X$ appel{\'e} espace des chemins d'ex{\'e}cution (non-constants),
de deux applications continues $s:\P X\rightarrow X^0$ et $t:\P X\rightarrow X^0$ ($X^0$ {\'e}tant muni
de la topologie discr{\`e}te) et d'une application continue
$*:\{(x,y)\in \P X\p\P X,t(x)=s(y)\}\rightarrow \P X$ satisfaisant les axiomes $s(x*y)=s(x)$,
$t(x*y)=t(y)$ et enfin
$x*(y*z)=(x*y)*z$ pour tout $x,y,z\in \P X$. Un morphisme de flots $f$ de $X$ vers $Y$ est une application
continue de $X^0\sqcup \P X$ vers $Y^0\sqcup \P Y$ tels que $f(X^0)\subset Y^0$, $f(\P X)\subset \P Y$,
$s(f(x))=f(s(x))$, $t(f(x))=f(t(x))$, $f(x*y)=f(x)*f(y)$. La cat{\'e}gorie correspondante est not{\'e}e
$\dtop$.
\end{definition}

Si $Z$ est un espace topologique, le flot $Glob(Z)$ est d{\'e}fini comme suit : $\P Glob(Z)=Z$,
$Glob(Z)^0=\{0,1\}$, et enfin $s=0$ et $t=1$ (il n'y a pas de chemins d'ex{\'e}cution composables).

\begin{definition}\cite{flow}
Deux morphismes de flots $f$ et $g$ de $X$ dans $Y$ sont S-homotopes s'il existe
une application continue $H:X\p [0,1]\rightarrow Y$ telle que
\begin{enumerate}
\item $H(-,u)$ est un morphisme de flots de $X$ dans $Y$ pour tout $u\in [0,1]$
\item $H(-,0)=f$ et $H(-,1)=g$
\end{enumerate}
On {\'e}crit $f\sim_S g$.
\end{definition}

On remarque que deux morphismes de flots S-homotopes co{\"\i}ncident sur le
$0$-squelette.

\begin{definition}\cite{flow}
Deux flots $X$ et $Y$ sont S-homotopes s'il existe un morphisme de
flots
$f:X\rightarrow Y$ et un morphisme de flots  $g:Y\rightarrow X$ tels que
$f\circ g\sim_{S} Id_Y$ et $g\circ f\sim_{S} Id_X$. \end{definition}

\begin{definition}\cite{flow}
Soit $X$ un flot et soit $Y$ un sous-ensemble de $X^0$.
La restriction $X\!\restriction_Y$
de $X$ {\`a} $Y$ est le flot d{\'e}fini par
\begin{enumerate}
\item le $0$-squelette de $X\!\restriction_Y$ est $Y$
\item $\P (X\!\restriction_Y)=\bigsqcup_{(\alpha,\beta)\in Y\p Y} \P_{\alpha,\beta}X$
 avec une d{\'e}finition {\'e}vidente des source, but et loi de composition.
\end{enumerate}
\end{definition}

Si $X$ est un flot, l'ensemble
$\P_\alpha^- X$ (resp.  $\P_\alpha^+ X$) est, dans la d{\'e}finition qui suit, l'ensemble des $\gamma\in
\P X$ tels que $s(\gamma)=\alpha$ (resp. $t(\gamma)=\alpha$) quotient{\'e}
par la relation d'{\'e}quivalence identifiant $\gamma$ avec
$\gamma*\gamma'$ (resp. $\gamma'*\gamma$ et $\gamma$). C'est-{\`a}-dire
que $\P_\alpha^- X$ (resp. $\P_\alpha^+ X$) est l'ensemble des classes
d'{\'e}quivalence de chemins d'ex{\'e}cution non-constants commen{\c c}ant (resp. terminant) de
la m{\^e}me fa{\c c}on en $\alpha$.

\begin{definition}\cite{flow}
Un morphisme de flots $f:X\rightarrow Y$ est une T-dihomotopie quand
\begin{enumerate}
\item Le morphisme de flots $X\rightarrow X\!\restriction_{f(X^0)}$ est un
isomorphisme de flots
\item pour $\alpha\in Y^0\backslash f(X^0)$, $\P^-_{\alpha}Y$ et $\P^+_{\alpha}Y$
sont des singletons
\item pour tout $\gamma\in Y\backslash f(X)$, il existe $u,v\in Y$ tels que $u*\gamma*v$
soit dans l'image $f(X)$ (avec les conventions $s(x)*x=x*t(x)=x$).
\end{enumerate}
On dit alors que les flots $X$ et $Y$ sont T-dihomotopes.
\end{definition}

\section{Comparaison entre CW-complexe globulaire et flot}

\begin{theoreme}\cite{flow}
Il existe un et un seul foncteur $cat$ (appel{\'e} r{\'e}alisation
cat{\'e}gori\-que) de $\gltop$ dans $\dtop$ tel que $cat(X^0)=cat(X)^0$,
$cat(Glob^{top}(Z))=Glob(Z)$, et $cat(X)=\limind_Z Glob(Z)$ o{\`u} $Z$
parcourt la d{\'e}composition cellulaire du CW-comple\-xe globulaire
$X$. \end{theoreme}

Si $X$ est un espace discret, on pose naturellement $cat(X):=X$. Si
$Z$ est un compact, on pose naturellement
$cat(Glob^{top}(Z)):=Glob(Z)$. Puis on d{\'e}finit le plongement sur les
objets en posant $cat(X):=\limind_Z Glob(Z)$ o{\`u} $Z$ parcourt la
d{\'e}composition cellulaire de $X$. Le plongement sur les morphismes est
obtenu comme suit. On construit d'abord explicitement pour tout
morphisme de CW-complexes globulaires $f:Glob^{top}(Z)\rightarrow U$ le
morphisme de flots $cat(f):Glob(Z)\rightarrow cat(U)$. Puis on pose
pour $f:X\rightarrow U$ un morphisme de CW-complexes globulaires
quelconque $cat(f):=\limind_Z cat(f\!\restriction_{Glob^{top}(Z)})$ o{\`u} $Z$
parcourt la d{\'e}composition cellulaire de $X$. On a donc par
construction la bijection $\dtop(cat(X),cat(U))\iso \limproj_n
\dtop(cat(X^1_n),cat(U))$.

\begin{theoreme}\label{sur}
\cite{flow} Soit $\tgltop(X,U)$ (resp. $\tdtop(cat(X),cat(U))$)
l'espace des morphismes de CW-complexes globulaires de $X$ vers $U$
(resp. de flots de $cat(X)$ vers $cat(U)$) muni de la Kelleyfication
de la topologie induite par celle de $\ttop(X,U)$
(resp. $\ttop(cat(X),cat(U))$). Alors $cat$ induit une homotopie
faible \[cat_*:\tgltop(X,U)\rightarrow
\tdtop(cat(X)\linebreak[0],\linebreak[0]cat(U)).\]  De plus il existe
une application continue $r$ de \[\tdtop(cat(X),cat(U))\] dans
$\tgltop(X,U)$ telle que $cat_*\circ r=Id_{\tdtop(cat(X),cat(U))}$.
\end{theoreme}

Consid{\'e}rons d'abord le cas $X=Glob^{top}(Z)$. La construction de $r$ se
fait alors explicitement et on d{\'e}montre que c'est un inverse
homotopique de $cat_*$. Cela est essentiellement d{\^u} au fait que par un
point donn{\'e} de $Glob^{top}(Z)$ autre que les deux points du $0$-squelette, il ne passe
qu'un et un seul chemin d'ex{\'e}cution non-constant ({\`a} reparam{\'e}trisation pr{\`e}s).

On obtient alors $r$ par it{\'e}ration sur la d{\'e}composition cellulaire de
$X$. De plus avec
\cite{0259.55004}
ou \cite{MR2001d:55012}, on en d{\'e}duit que $\holimproj_n
\tgltop(X^1_n\linebreak[0],\linebreak[0]U)$ est faiblement
homotope {\`a} $\holimproj_n \tdtop(cat(X^1_n),cat(U))$.

On montre ensuite que les deux tours d'espaces donnant les deux
limites projectives homotopiques sont des tours fibrantes, {\`a} savoir
que toutes les applications apparaissant dans ces tours sont des
fibrations d'espaces topologiques. Cela est essentiellement d{\^u} au
fait que les inclusions $S^{n-1}\hookrightarrow D^n$ et donc
$Glob^{top}(S^{n-1})\hookrightarrow Glob^{top}(D^n)$ sont des cofibrations
\cite{MR80b:55001}. On en d{\'e}duit donc, encore gr{\^a}ce {\`a}
\cite{0259.55004} ou \cite{MR2001d:55012}, les {\'e}quivalences
d'homotopie faible $\limproj_n \tgltop(X^1_n,U)\iso
\holimproj_n \tgltop(X^1_n,U)$ et $\limproj_n
\tdtop(cat(\linebreak[0]X^1_n\linebreak[0])\linebreak[0],\linebreak[0]cat(\linebreak[0]U\linebreak[0])\linebreak[0]\linebreak[0])
\iso \holimproj_n \tdtop(cat(X^1_n),cat(\linebreak[0]U\linebreak[0])\linebreak[0])$. Le r{\'e}sultat
d{\'e}coule alors des hom{\'e}omorphismes \[\tgltop(X\linebreak[0],\linebreak[0]U\linebreak[0])\linebreak[0]\iso\linebreak[0]\limproj_n \tgltop(\linebreak[0]X^1_n,U)\]
et \[\tdtop(cat(X),cat(U))\iso \limproj_n \tdtop(cat(X^1_n)\linebreak[0],\linebreak[0]cat(U)).\]

\begin{theoreme}\cite{flow} Deux CW-complexes globulaires $X$ et $Y$ sont S-homotopes
si et seulement si les flots $cat(X)$ et $cat(Y)$ sont S-homotopes.
\end{theoreme}

Les espaces topologiques $\tgltop(X,U)$ et $\tdtop(cat(X),cat(U))$
sont faiblement homotopes, et ont donc les m{\^e}mes composantes connexes
par arc. Comme on travaille dans une cat{\'e}gorie d'espaces topologiques
qui est cart{\'e}siennement ferm{\'e}e, une composante connexe par arc de
$\tgltop(X,U)$
(resp. $\tdtop(\linebreak[0]cat(\linebreak[0]X)\linebreak[0],\linebreak[0]cat(U))$)
est une classe de S-homotopie de morphismes de CW-complexes
globulaires (resp. de flots). Le r{\'e}sultat d{\'e}coule alors de la
surjectivit{\'e} de l'application $cat_*$.

Enfin on a aussi :

\begin{theoreme}\cite{flow} Deux CW-complexes globulaires $X$ et $Y$ sont T-homotopes
si et seulement si les flots $cat(X)$ et $cat(Y)$ sont T-homotopes.
\end{theoreme}

\providecommand{\bysame}{\leavevmode\hbox to3em{\hrulefill}\thinspace}
\providecommand{\MR}{\relax\ifhmode\unskip\space\fi MR }
\providecommand{\MRhref}[2]{%
  \href{http://www.ams.org/mathscinet-getitem?mr=#1}{#2}
}
\providecommand{\href}[2]{#2}

\end{document}